\documentclass{amsart}
\usepackage{amsmath} 
\usepackage{amssymb}
\usepackage{mathrsfs}

\newtheorem{theorem}{Theorem}[section] 
\newtheorem{claim}[theorem]{Claim}

\newtheorem{conclusion}[theorem]{Conclusion}

\theoremstyle{definition}
\newtheorem{definition}[theorem]{Definition}

\theoremstyle{remark}
\newtheorem{remark}[theorem]{Remark}

\newcommand{\ch}{{\rm ch}}

\newcommand{\Ord}{{\rm Ord}}

\newcommand{\INC}{{\rm INC}}

\newcommand{\Inc}{{\rm Inc}}

\newcommand{\chr}{{\rm chr}}

\newcommand{\edge}{{\rm edge}}

\newcommand{\Dom}{{\rm Dom}}
\newcommand{\Rang}{{\rm Rang}}

\newcommand{\rest}{{\restriction}}

\newcommand{\then}{{\underline{then}}}
\newcommand{\when}{{\underline{when}}}

\newcommand{\where}{{\underline{where}}}

\newcommand{\mn}{{\medskip\noindent}}
\newcommand{\sn}{{\smallskip\noindent}}

\newcommand{\cA}{{\mathscr A}}

\newcommand{\cB}{{\mathscr B}}

\newcommand{\cH}{{\mathscr H}}

\newcommand{\cF}{{\mathscr F}}

\newcommand{\bbZ}{{\mathbb Z}}

\newcommand{\gt}{{\mathfrak t}}

\newcommand{\pr}{{\rm pr}}

\newcount\skewfactor
\def\mathunderaccent#1#2 {\let\theaccent#1\skewfactor#2
\mathpalette\putaccentunder}
\def\putaccentunder#1#2{\oalign{$#1#2$\crcr\hidewidth
\vbox to.2ex{\hbox{$#1\skew\skewfactor\theaccent{}$}\vss}\hidewidth}}

\newenvironment{PROOF}[2][\proofname.]
   {\begin{proof}[#1]\renewcommand{\qedsymbol}{\textsquare$_{\rm #2}$}}
   {\end{proof}}

\begin{document}

\title {More on compactness of chromatic number}
\author {Saharon Shelah}
\address{Einstein Institute of Mathematics\\
Edmond J. Safra Campus, Givat Ram\\
The Hebrew University of Jerusalem\\
Jerusalem, 91904, Israel\\
 and \\
 Department of Mathematics\\
 Hill Center - Busch Campus \\ 
 Rutgers, The State University of New Jersey \\
 110 Frelinghuysen Road \\
 Piscataway, NJ 08854-8019 USA}
\email{shelah@math.huji.ac.il}
\urladdr{http://shelah.logic.at}
\thanks{The author would like to thank the Israel Science Foundation
  for partial support of this research (Grant No. 1053/11).
The author thanks Alice Leonhardt for the beautiful typing.
First typed July 3, 2011. Publication 1018.}
 
 






\date{February 4, 2013}

\begin{abstract}
We prove that for any regular $\kappa$ and $\mu > \kappa$ below the first
fix point $(\lambda = \aleph_\lambda)$ above $\kappa$, there is a graph
with chromatic number $> \kappa$, and $\mu^\kappa$ nodes 
but every subgraph of cardinality $< \mu$ has chromatic number $\le \kappa$.
\end{abstract}

\maketitle
\numberwithin{equation}{section}
\setcounter{section}{-1}
\newpage

\section {Introduction}

This continues \cite{Sh:1006} but does not rely on it.

In \cite{Sh:1006} we prove that if there is $\cF \subseteq
{}^\kappa\Ord$ of cardinality $\mu,\lambda$-free not free then we can
get a failure of $\lambda$-compactness for the chromatic number being
$\kappa$.  This gives (using \cite[Ch.II]{Sh:g}) that if 
$\mu$ is strong limit singular of
cofinality $\kappa$ and $2^\mu > \mu^+$ then we get the above for
$\lambda = \mu^+$ (and more).

Now by \cite{Sh:161} having a $\lambda$-free not free Abelian group of
cardinality $\lambda$ is characterized combinatorically, in particular by
freeness for the existence of transversals.
By Magidor-Shelah \cite{MgSh:204}, for arbitrarily large $\lambda <$
first $\mu = \aleph_\mu$, there is a $\lambda$-free, not free Abelian
group (of cardinality $\lambda$), see history there.  We prove, derived
from such incompactness examples for incompactness of being
$\kappa$-chromatic, in particular answering a second problem of
Magidor: $\aleph_\omega$-compactness fails for being
$\aleph_0$-chromatic; moreover this holds for $\mu <$ first fix point
of the $\aleph$'s.

We intend to continue in \cite{Sh:F1296}.

Another problem on incompactness is about the existence of
$\lambda$-free Abelian groups $G$ which with no non-trivial
homomorphism to $\bbZ$, in \cite{Sh:883}, for $\lambda = \aleph_n$
using $n-BB$.  
In \cite{Sh:898} we get more $\lambda$'s, almost in ZFC by 1-BB
(black box).  This proof suffices here (but not in ZFC).
This is continued in \cite{Sh:F1200}, but presently not connected).

\begin{definition}
\label{y8}
1) We say ``we have $(\mu,\lambda)$-incompactness for the $(< \chi)$-chromatic
 number" or $\INC_{\chr}(\mu,\lambda,< \chi)$
 \when \, there is an increasing continuous sequence $\langle
 G_i:i \le \lambda\rangle$ of graphs each with $\le \mu$ nodes, $G_i$ an
 induced subgraph of $G_\lambda$ with $\ch(G_\lambda) \ge \chi$ but $i
 < \lambda \Rightarrow \ch(G_i) < \chi$.

\noindent
2) Replacing (in part (1)) $\chi$ by $\bar\chi = (< \chi_0,\chi_1)$
means $\ch(G_\lambda)) \ge \chi_1$ and $i < \lambda \rightarrow \ch(G_i) <
\chi_0$; similarly in parts 3),4) below.

\noindent
3) We say we have incompactness for length $\lambda$ for $(< \chi)$-chromatic
(or $\bar\chi$-chromatic) number \when \, we fail to have
$(\mu,\lambda)$-compactness for $(< \chi)$-chromatic (or
   $\bar\chi$-chromatic) number for some $\mu$.

\noindent
4) We say we have $[\mu,\lambda]$-incompactness for $(<
   \chi)$-chromatic number or $\INC_{\chr}[\mu,\lambda,< \chi]$ 
\when \, there is a graph $G$ with $\mu$ nodes, 
$\ch(G) \ge \chi$ but $G^1 \subseteq G \wedge |G^1| < \lambda \Rightarrow
   \ch(G^1) < \chi$.

\noindent
5) Let $\INC^+_{\chr}(\mu,\lambda,< \chi)$ be as in part (1) but we
add that there is a partition $\langle A_{1,\varepsilon}:\varepsilon <
\kappa\rangle$ of the set of nodes of $G_i$ such that
$c \ell(G_i \rest A_{i,\varepsilon})$, the colouring number of $G_i
\rest A_{i,\varepsilon}$ is $<  \chi$ for $i < \lambda$, see below.

\noindent
6) Let $\INC^+_{\chr}[\mu,\lambda,< \chi]$ be as in part (4) but we
add: if $G^1 \subseteq G$ and $|G^1| < \lambda$ \then \, there is a
partition $\langle A_\varepsilon:\varepsilon < \varepsilon_*\rangle$
of the nodes of $G^1$ to $\varepsilon_* < \chi$ sets such that $\varepsilon <
\varepsilon_* \Rightarrow c \ell(G^1 \rest A_\varepsilon) < \chi$.

\noindent
7) If $\chi = \kappa^+$ we may write $\kappa$ instead of ``$< \chi$".

\noindent
8) Let $\INC(\lambda,< \chi)$ means $\INC(\lambda,\lambda,< \chi)$,
   and similarly in the other cases.
\end{definition}
\newpage

\section {A sufficient criterion and relations to transversals}

\begin{definition}
\label{r2}
1) Let $\Inc[\mu,\lambda,\kappa]$ mean that we can find $\bold a = 
(\cA,\bar R)$ witnessing it which means that:
\mn
\begin{enumerate}
\item[$(a)$]  $|\cA| = \mu$
\sn
\item[$(b)$]  $\bar R = \langle R_\varepsilon:\varepsilon
  <\kappa\rangle$
\sn
\item[$(c)$]  $R_\varepsilon$ is a two-place relation on $\cA$, so we
  may write $\nu R_\varepsilon \eta$
\sn
\item[$(d)$]  $\cA$ is not free (for $\bold a$), see $(*)_1$ below or just
 not strongly free, see $(*)_2$ below
\sn
\item[$(e)$]  $\bold a = (\cA,\bar R)$ is $\lambda$-free which means
  $\cB \subseteq \cA \wedge |\cB| < \lambda \Rightarrow \cB$ is $\bold a$-free
\end{enumerate}
\mn
\where
\mn
\begin{enumerate}
\item[$(*)_1$]  if $\cB \subseteq \cA$ \then \,
  \, $\cB$ is $\bold a$-free means that there is a witness $(h,<_*)$
  which means
\sn
\item[${{}}$]  $(\alpha) \quad <_*$ a well ordering of $\cB$
\sn
\item[${{}}$]  $(\beta) \quad h$ is a function from $\cB$ to $\kappa$
\sn
\item[${{}}$]  $(\gamma) \quad$ if $h(\eta) = h(\nu)$ and $\nu R_\zeta
  \eta$ for some $\zeta$ \then \, $\nu <_* \eta$ (so really only

\hskip25pt  $<_* \rest \{\eta \in \cB:
h(\eta) = \varepsilon\}$ for $\varepsilon < \kappa$ count); so it is
reasonable 

\hskip25pt to assume each $R_\varepsilon$ is irreflexive
\sn
\item[${{}}$]  $(\delta) \quad$ for any $\eta \in \cB$ the
set\footnote{exp stands for exceptional}
  $\exp(\eta,h,<_*)$ has cardinality $< \kappa$
where 

\hskip25pt (recall that $\cB = \Dom(h))$
 
\hskip25pt $\bullet \quad \exp(\eta,h,<_*) = \exp(\eta,h,<_*,\bold a) 
= \{\zeta < \kappa$: there is $\nu <_* \eta$ 

\hskip40pt  such that $\nu R_\zeta \eta$ and $h(\nu) = h(\eta)\}$
\sn
\item[$(*)_2$]  if $\cB \subseteq \cA$ \then \,
$\cB$ is strongly $\bold a$-free means that for every well ordering
  $<_*$ of $\cB$ there is a function $h:\cB \rightarrow \kappa$ such
  that $(h,<_* \rest \cB)$ witness $\cB$ is $\bold a$-free
\sn
\item[$(*)_3$]  if $\cB \subseteq \cA$ \then \, $\cB$ is weakly free means
that there is a witness $h$ which means
\sn
\item[${{}}$]  $(\alpha) \quad h$ is a function from $\cB$ to $\kappa$
\sn
\item[${{}}$]  $(\beta) \quad$ for every $\eta \in \cB$ the set
$\exp(\eta,h)$ has cardinality $< \kappa$ where
\sn
\item[${{}}$]  $\qquad \bullet \quad \exp(\eta,h) = \exp(\eta,h,\bold
a) = \{\zeta < \kappa$: there is $\nu \in \cB$ such that $\nu R_\zeta
\eta$

\hskip40pt  and $h(\nu) = h(\eta)\}$.
\end{enumerate}
\mn
2) Let $\Inc(\mu,\lambda,\kappa)$ mean that we can find
$(\cA,\bar{\cA},\bar R)$ witnessing it which means that:
\mn
\begin{enumerate}
\item[$(a)-(d)$]  as above
\sn
\item[$(e)'$]  $\bar{\cA} = \langle \cA_\alpha:\alpha < \lambda\rangle$
  is an increasing sequence with union $\cA$ such that for each
  $\alpha < \lambda$ the set $\cA_\alpha$ is free (i.e. for $(\cA,\bar R)$).
\end{enumerate}
\end{definition}

\begin{claim}
\label{r6}
We have $\INC_{\chr}(\mu,\lambda,\kappa)$ or
$\INC_{\chr}[\mu,\lambda,\kappa]$, see Definition \ref{y8}(4)
 \when \,:
\mn
\begin{enumerate}
\item[$\boxplus$]  $(a) \quad \Inc(\chi,\lambda,\kappa)$ or
  $\Inc[\chi,\lambda,\kappa]$ respectively
\sn
\item[${{}}$]  $(b) \quad \chi \le \mu = \mu^\kappa$.
\end{enumerate}
\end{claim}

\begin{PROOF}{\ref{r6}}
Fix $\bold a = (\cA,\bar{\cA},\bar R)$ or 
$\bold a = (\cA,\bar R)$ witnessing $\Inc(\mu,\lambda,\kappa)$ 
or $\Inc[\mu,\lambda,\kappa]$ respectively.  Now we define $\tau_{\cA}$ as the
vocabulary $\{P_\eta:\eta \in \cA\} \cup \{F_\varepsilon:\varepsilon <
\kappa\}$ where $P_\eta$ is a unary predicate, $F_\varepsilon$ a unary
function (but it may be interpreted as a partial function).

We further let $K_{\bold a}$ be the class of structures $M$ such that:
\mn
\begin{enumerate}
\item[$\boxplus_1$]  $(a) \quad M = (|M|,F^M_\varepsilon,
P^M_\eta)_{\varepsilon < \kappa,\eta \in \cA}$
\sn
\item[${{}}$]  $(b) \quad \langle P^M_\eta:\eta \in \cA\rangle$ is a
partition of $|M|$, so for $a \in M$ let $\eta[a]$ 

\hskip25pt $= \eta^M_a$ be the unique $\eta \in \cA$ such that $a \in P^M_\eta$
\sn
\item[${{}}$]  $(c) \quad$ if $a_\ell \in P^M_{\eta_\ell}$ for
  $\ell=1,2$ and $F^M_\zeta(a_2) = a_1$ then

\hskip25pt $\eta_1 R_\zeta \eta_2$.
\end{enumerate}
\mn
Let $K^*_{\bold a}$ be the class of $M$ such that
\mn
\begin{enumerate}
\item[$\boxplus_2$]  $(a) \quad M \in K_{\bold a}$
\sn
\item[${{}}$]  $(b) \quad \|M\| = \mu$
\sn
\item[${{}}$]  $(c) \quad$ if $\eta \in \cA,u \subseteq \kappa$ and
  for $\zeta \in u$ we have 
$\nu_\zeta \in \cA,\nu_\zeta R_\zeta \eta$ and $a_\zeta \in 
P^M_{\nu_\zeta}$ 

\hskip25pt   \then \, for some $a \in P^M_\eta$ we have 
$\zeta \in u \Rightarrow F^M_\zeta(a) = a_\zeta$ 

\hskip25pt and $\zeta \in \kappa \backslash u \Rightarrow 
F^M_\zeta(a)$ not defined.
\end{enumerate}
\mn
Clearly
\mn
\begin{enumerate}
\item[$\boxplus_3$]   there is $M \in K^*_{\bold a}$.
\end{enumerate}
\mn
[Why?  Obvious as we are assuming $|\cA| = \chi \le \mu = \mu^\kappa$.]
\mn
\begin{enumerate}
\item[$\boxplus_4$]  for $M \in K_{\bold a}$ let $G_M$ be the graph with:
\sn
\begin{enumerate}
\item[$\bullet$]  set of nodes $|M|$
\sn
\item[$\bullet$]  set of edges $\{\{a,F^M_\varepsilon(a)\}:a \in
  |M|,\varepsilon < \kappa$ when $F^M_\varepsilon(a)$ is defined$\}$.
\end{enumerate}
\end{enumerate}
\mn
We shall show that the graph $G_M$ is as required in Definition
\ref{y8}(1) or \ref{y8}(4) (recalling $\kappa^+$ here stands for
$\chi$ there, see \ref{y8}(7).  Clearly $G_M$ is a graph with $\mu$
nodes so recalling Definition \ref{r2}(2) or \ref{r2}(1) it suffices
to prove $\boxplus_5$ and $\boxplus_7$ below.
\mn
\begin{enumerate}
\item[$\boxplus_5$]  if $\cB \subseteq \cA$ 
is free, and $M \in K_{\bold a}$ \then \, $G_{M,\cB} := G_M \rest 
(\cup\{P^M_\eta:\eta \in \cB\})$
has chromatic number $\le \kappa$.
\end{enumerate}
\mn
[Why?  Let the pair $(h,<_*)$ witness that $\cB$ is free 
(for $\bold a = (\cA,\bar R)$, see \ref{r2}(1)$(*)_1$)
so $h:\cB \rightarrow \kappa$
and let $\cB_\varepsilon = \{\eta \in \cB:h(\eta) = 
\varepsilon\}$ for $\varepsilon < \kappa$.  

Clearly
\mn
\begin{enumerate}
\item[$\boxplus_{5.1}$]  it suffices for each $\varepsilon < \kappa$
  to prove that $G_{M,\cB_\varepsilon}$ has chromatic number $\le
  \kappa$.
\end{enumerate}
\mn
Let $\langle \eta_\alpha:\alpha < \alpha(*)\rangle$ list $\cB$ in
 $<_*$-increasing order.
We define $\bold c_\varepsilon:G_{M,\cB_\varepsilon} \rightarrow
\kappa$ by defining a colouring $\bold c_{\varepsilon,\alpha}:
G_{M,\{\eta_\beta:\beta < \alpha\} \cap \cB_\varepsilon} \rightarrow
\kappa$ by induction on $\alpha \le \alpha(*)$
such that $\bold c_{\varepsilon,\alpha}$ 
is increasing continuous with $\alpha$.  
For $\alpha = 0$, let $\bold c_{\varepsilon,\alpha} = \emptyset$, and
for $\alpha$ limit take union.  If $\alpha = \beta +1$ and $\eta_\beta
\notin \cB_\varepsilon$ then we let $\bold c_\alpha = \bold c_\beta$.

Lastly, assume $\alpha = \beta +1,\eta_\beta \in \cB_\varepsilon$ \then
\,note that the set $u_{\varepsilon,\beta} = \{\zeta < \kappa$: 
there is $\nu <_* \eta_\beta$ such that $\nu \in
\cB_\varepsilon$ and $\nu R_\zeta \eta\}$ has cardinality 
$< \kappa$ because the pair $(<_*,h)$ witness ``$\cB$ is free".  
Hence, recalling $M \in K_{\bold a}$, for each 
$a \in P^M_{\eta_\beta}$, the set 
$u_{\varepsilon,\beta,a} := \{\zeta < \kappa_\varepsilon:
F^M_\zeta(a) \in \{P^M_\nu:\nu <_* \eta_\beta$
and $\nu \in \cB_\varepsilon\}\}$ is $\subseteq u_{\varepsilon,\beta}$
hence has cardinality $\le |u_{\varepsilon,\beta}| < \kappa$.  But by
$(*)_1(\gamma)$ of \ref{r2} and the definition of $K_{\bold a},
A_a := \{b \in G_{M,\{\eta_\gamma:\gamma < \beta\} \cap 
\cB_\varepsilon}:\{b,a\}$ is an edge of $G_M\}$ is $\subseteq
\{F^M_\zeta(a):\zeta \in u_{\varepsilon,\beta,a}\}$ 
hence the set $A_a$ has cardinality $\le |u_{\varepsilon,\beta,a}| <\kappa$.
So define $\bold c_{\varepsilon,\alpha}$ extending $\bold
c_{\varepsilon,\beta}$ by, for $a \in P^M_{\eta_\beta}$ letting
 $\bold c_{\varepsilon,\alpha}(a) = \min(\kappa \backslash 
\{\bold c_{\varepsilon,\beta}(b):b \in P^M_\nu$ for some $\nu <_*
\eta_\beta$ from $\cB_\varepsilon$ and $\{b,a\}$ is an edge of
$G_M\})$.   Recalling there is no edge
$\subseteq P_{\eta_\beta}$ this is a colouring.

So we can carry the induction. So indeed $\boxplus_5$ holds.]
\mn
\begin{enumerate}
\item[$\boxplus_6$]  if $\cB \subseteq \cA$ is free and $M \in
  K_{\bold a}$ then $G_{M,\cB}$ is the union of $\le \kappa$ sets each
with colouring number $\le \kappa$ hence also chromatic number $\le \kappa$.
\end{enumerate}
\mn
[Why?  By the proof of $\boxplus_5$.]
\mn
\begin{enumerate}
\item[$\boxplus_7$]  $\chr(G_M) > \kappa$ if $M \in K^*_{\cA}$.
\end{enumerate}
\mn
Why?  Toward contradiction assume $\bold c:G_M \rightarrow \kappa$ is
a colouring and let $<_*$ be a well ordering of $\cA$.  
For each $\eta \in \cA$ and $\varepsilon,\zeta < \kappa$ let
$\Lambda_{\eta,\varepsilon,\zeta} = \{\nu:\nu \in \cA,\nu <_* \eta,
\nu R_\zeta \eta$ and $\varepsilon \in \cH_\nu\}$ where for $\nu \in
\cA$ we define $\cH_\nu = \{\varepsilon$: for some $a \in
P^M_\nu$ we have $\bold c(a) = \varepsilon\}$.
\medskip

\noindent
\underline{Case 1}:  There is $\eta \in \cA$ such that $(\forall
\varepsilon \in \cH_\eta)(\exists^\kappa \zeta < \kappa)
[\Lambda_{\eta,\varepsilon,\zeta} \ne \emptyset]$. 

So we can find a one-to-one function $g:\cH_\eta \rightarrow \kappa$
such that $\Lambda_{\eta,\varepsilon,g(\varepsilon)} \ne \emptyset$
for every $\varepsilon \in \cH_\eta \subseteq \kappa$.
For each $\varepsilon \in \cH_\eta \subseteq \kappa$ 
choose $\nu_\varepsilon \in
\Lambda_{\eta,\varepsilon,g(\varepsilon)}$; 
possible as $\Lambda_{\eta,\varepsilon,g(\varepsilon)} \ne \emptyset$ 
by the choice of the function $g$.  By the definition of ``$\nu_\varepsilon \in
\Lambda_{\eta,\varepsilon,g(\varepsilon)}$" there is $a_\varepsilon \in
P^M_{\nu_\varepsilon}$ such that $\bold c(\nu_\varepsilon) =
\varepsilon$ and $\nu R_\zeta \eta$ holds.  
So as $M \in K^*_{\bold a}$ there is $a \in P^M_\eta$ 
such that $\varepsilon \in \cH_\eta \subseteq \kappa \Rightarrow
F^M_{g(\varepsilon)}(a) = a_\varepsilon$, but then $\{a,a_\varepsilon\} \in
\edge(G_M)$ hence $\bold c(a) \ne \bold c(a_\varepsilon) =
\varepsilon$ for every $\varepsilon \in \cH_\eta \subseteq 
\kappa$, contradiction to the definition of $\cH_\eta$. 
\medskip

\noindent
\underline{Case 2}:  Not Case 1

So for every $\eta \in \cA$ there is $\varepsilon \in \cH_\eta \subseteq
\kappa$ such that there are
$< \kappa$ ordinals $\zeta < \kappa$ such that
  $\Lambda_{\eta,\varepsilon,\zeta} \ne \emptyset$.  This means that
there is $h:\cA \rightarrow \kappa$ such that:
\mn
\begin{enumerate}
\item[$\bullet_1$]  $\eta \in \cA \Rightarrow h(\eta) \in \cH_\eta$
and
\sn
\item[$\bullet_2$]  $\eta \in \cA \Rightarrow \kappa > 
|\{\zeta < \kappa:\Lambda_{\eta,h(\eta),\zeta} \ne \emptyset\}|$.
\end{enumerate}
\mn
This implies that:
\mn
\begin{enumerate}
\item[$\bullet_3$]  $\eta \in \cA \Rightarrow \kappa >
  |\exp(\eta,h,\bold a,<_*)|$ 
\end{enumerate}
\mn
because
\mn
\begin{enumerate}
\item[$\bullet_4$]  if $\eta \in \cA$ and $\varepsilon = h(\eta)$
  \then \, $\exp(\eta,h,<_*,\bold a) \subseteq \{\zeta <
  \kappa:\Lambda_{\eta,\varepsilon,\zeta} \ne \emptyset\}$.
\end{enumerate}
\mn
[Why?  As $h:\cA \rightarrow \kappa$ and if $\zeta \in
  \exp(\eta,h,<_*,\bold a)$ 
let $\nu$ exemplify this, that is, $\nu <_* \eta,\nu R_\zeta \eta$ 
and $h(\nu) = h(\eta) = \varepsilon$ and recall $h(\nu) = \varepsilon$
implies $\varepsilon \in \cH_\nu$.  But this means
  that $\nu \in \Lambda_{\eta,\varepsilon,\zeta}$ hence
  $\Lambda_{\eta,\varepsilon,\eta} \ne \emptyset$ as required.]

As $<_*$ was any well ordering of $\cA$, this means, 
see \ref{r2}$(*)_2$, that $\cA$ is strongly free, contradiction to
\ref{r2}(d). 
\end{PROOF}

\begin{claim}
\label{r7}
1) We have $\Inc(\mu,\lambda,\kappa)$ \when \,
\mn
\begin{enumerate}
\item[$(*)$]  for some $\cF$ and natural number $\bold k > 0$ we have
\sn
\begin{enumerate}
\item[$(a)$]  $\cF \subseteq {}^\kappa \mu$ has\footnote{may add
$(b)$ \quad [normality] if $\eta,\nu \in \cF$ and $i,j < \kappa$ and
$\eta(i) = \nu(j)$ \then \, $i=j$ hence for every $\eta,\nu \in
  \cF$ we have $\eta = \nu \Leftrightarrow \Rang(\eta) = \Rang(\nu)$.}
cardinality $\mu$
\sn
\item[$(b)$]  $\cF$ is not free where
\sn
\item[${{}}$]  $\bullet \quad \cF' \subseteq \cF$ is free \when
\sn
\item[${{}}$]  $\bullet \quad$ there is a sequence 
$\langle \cF'_i:i < \kappa\rangle$ such
  that $\cF' = \cup\{\cF'_i:i < \kappa\}$ and for each $i,\cF'_i$ has
  a transversal which means that $\{\Rang(\eta):\eta \in \cF'_i\}$ has a
 transversal (= one-to-one choice function)
\sn
\item[$(c)$]  $\cF$ is the increasing union of $\langle
 \cF_\alpha:\alpha < \lambda\rangle$ such that each $\cF_\alpha$ is free.
\end{enumerate}
\end{enumerate}
\mn
2) We have $\Inc[\mu,\lambda,\kappa]$ \when \,
\mn
\begin{enumerate}
\item[$(*)$]  as above but replacing clause (c) by:
\sn
\item[${{}}$]  $(c)' \quad$ every $\cF' \subseteq \cF$ of cardinality
  $< \lambda$ has a transversal.
\end{enumerate}
\end{claim}

\begin{PROOF}{\ref{r7}}
1), 2) We define $\bold a$ by choosing (for our $\cF$):
\mn
\begin{enumerate}
\item[$\bullet$]  $\cA_{\bold a} = \cF$
\sn
\item[$\bullet$]  $<_{\cA}$ any well ordering of $\cF$; not part of
$\bold a$
\sn
\item[$\bullet$]  $R_\varepsilon$ is defined by: $f R_\varepsilon g$
  iff $f <_{\cA} g \wedge f(\varepsilon) = g(\varepsilon)$
\sn
\item[$\bullet$]  for part (1) let $\bar{\cA}$ be a sequence
  witnessing clause (c).
\end{enumerate}
\mn
So it suffices to prove $\Inc(\mu,\lambda,\kappa)$ or
$\Inc[\mu,\lambda,\kappa]$; hence it suffices to 
prove that $\bold a$ witness it.

Now in Definition \ref{r2}, clauses (a),(b),(c) are obvious.  For
clause (e), assume $\cF_2 \subseteq \cF$ is free in the sense of
\ref{r7}(1)(b), and we shall prove that $\cF_2$ is $\bold a$-free,
this suffices for clause (e).  By the assumption on $\cF_2$, clearly
$\cF_2$ is the union of $\langle \cF_{2,\zeta}:\zeta < \kappa\rangle,
\cF_{2,\zeta}$ has a transversal $\bold h_\zeta$.  Now we
define $h:\cF_2 \rightarrow \kappa$ by: $h(f) = \pr(\zeta,\varepsilon)$ where
$\zeta = \min\{\xi:f \in \cF_{2,\xi}\}$ and $\varepsilon$ is minimal
such that $\bold h_\zeta(\Rang(f)) = f(\varepsilon)$, now the pairs $(h,<_{\cA} \rest
\cF_2)$ witness that $\cF_2$ is free (for $\bold a$).

For clause (d) toward contradiction assume that 
$h:\cF \rightarrow \kappa$ and well ordering $<_*$ of $\cA$ witness $\cF$ is
free for $\bold a$, hence $\bar\cB = \langle \cB_\varepsilon:\varepsilon <
\kappa \rangle$ is a partition of $\cF$ when we let $\cB_\varepsilon = \{f \in
\cF:h(f) = \varepsilon\}$.

By Definition \ref{r2}, for each 
$\varepsilon < \kappa$ and $f \in \cB_\varepsilon$ the set
$u_f = \{\zeta < \kappa$: for some $g \in \cB_\varepsilon$ we have $g
R_\zeta f\}$ has cardinality $< \kappa$ and let $\zeta_f \in \kappa
\backslash u_f$.  For $\varepsilon,\zeta < \kappa$ let
$\cB_{\varepsilon,\zeta} = \{f \in \cB_\varepsilon:\zeta_f = \zeta\}$
so $\langle \cB_{\varepsilon,\zeta}:\varepsilon,\zeta < \kappa\rangle$
is a partition of $\cA$.  Now for each $\varepsilon,\zeta < \kappa$,
if $f \ne g \in \cB_{\varepsilon,\zeta}$ \then \, $f(\zeta) \ne
g(\zeta)$.  Why?  By symmetry we can assume $g <_{\cA} f$ now $\zeta =
\zeta_f \in \kappa \backslash u_f$, so $g$ cannot witness $\zeta \in u_f$.
So $\langle \cB_{\varepsilon,\zeta}:\varepsilon,\zeta <
\kappa\rangle$ contradicts clause (b) of the claim's assumption.
\end{PROOF}

\begin{conclusion}
\label{r8}
For any $\kappa$, for arbitrarily large regular $\lambda <
\min\{\theta:\theta = \aleph_\theta > \kappa\}$ we have
$\INC^+(\lambda,\lambda,\kappa)$. 
\end{conclusion}

\begin{remark}
\label{r12}
Note the variant of transversal we use for our purpose is equivalent
by \cite{Sh:161}.
\end{remark}

\begin{PROOF}{\ref{r8}}
If $\kappa = \aleph_0$, by Magidor-Shelah \cite{MgSh:204}.  Generally
similar.
\end{PROOF}

\begin{claim}
\label{r18}
If $\INC[\mu,\lambda,\kappa]$ or $\INC(\mu,\lambda,\kappa)$ \then
   \, $\Inc[\mu,\lambda,\kappa]$ or $\Inc(\mu,\lambda,\kappa)$ respectively.
\end{claim}

\begin{PROOF}{\ref{r18}}
As the two cases are similar we do the $\INC(\mu,\lambda,\kappa)$ case, so
let $G,\langle G_i:i < \lambda\rangle$ witness it.

Let $<_*$ be a well ordering of the set of nodes of $G$.  Define
$\bold a = (\cA,\bar{\cA},\bar R)$ by:
\mn
\begin{enumerate}
\item[$\bullet$]  $\cA$ is the set of nodes of $G$
\sn
\item[$\bullet$]  $\bar{\cA} = \langle \cA_i:i < \lambda \rangle$ with 
$\cA_i$ the set of nodes of $G_i$
\sn
\item[$\bullet$]  $R_\varepsilon = \{(\nu,\eta):\{\nu,\eta\}$ an edge
of $G$ and $\nu <_* \eta\}$.
\end{enumerate}
\mn
Now check, noting when checking,  that e.g. in $(*)_1$ of Definition
\ref{r2},  $\exp(\eta,\alpha,<_*)$ is equal to $\kappa$ or to
$\emptyset$ as $\bigwedge\limits_{\varepsilon} R_\varepsilon = R_0$.
\end{PROOF}
\newpage
\def\germ{\frak} \def\scr{\cal} \ifx\documentclass\undefinedcs
  \def\bf{\fam\bffam\tenbf}\def\rm{\fam0\tenrm}\fi 
  \def\defaultdefine#1#2{\expandafter\ifx\csname#1\endcsname\relax
  \expandafter\def\csname#1\endcsname{#2}\fi} \defaultdefine{Bbb}{\bf}
  \defaultdefine{frak}{\bf} \defaultdefine{=}{\B} 
  \defaultdefine{mathfrak}{\frak} \defaultdefine{mathbb}{\bf}
  \defaultdefine{mathcal}{\cal}
  \defaultdefine{beth}{BETH}\defaultdefine{cal}{\bf} \def\bbfI{{\Bbb I}}
  \def\mbox{\hbox} \def\text{\hbox} \def\om{\omega} \def\Cal#1{{\bf #1}}
  \def\pcf{pcf} \defaultdefine{cf}{cf} \defaultdefine{reals}{{\Bbb R}}
  \defaultdefine{real}{{\Bbb R}} \def\restriction{{|}} \def\club{CLUB}
  \def\w{\omega} \def\exist{\exists} \def\se{{\germ se}} \def\bb{{\bf b}}
  \def\equivalence{\equiv} \let\lt< \let\gt>
  \def\implies{\Rightarrow}\def\mathfrak{\bf}\def\germ{\frak} \def\scr{\cal}
  \ifx\documentclass\undefinedcs
  \def\bf{\fam\bffam\tenbf}\def\rm{\fam0\tenrm}\fi 
  \def\defaultdefine#1#2{\expandafter\ifx\csname#1\endcsname\relax
  \expandafter\def\csname#1\endcsname{#2}\fi} \defaultdefine{Bbb}{\bf}
  \defaultdefine{frak}{\bf} \defaultdefine{=}{\B} 
  \defaultdefine{mathfrak}{\frak} \defaultdefine{mathbb}{\bf}
  \defaultdefine{mathcal}{\cal}
  \defaultdefine{beth}{BETH}\defaultdefine{cal}{\bf} \def\bbfI{{\Bbb I}}
  \def\mbox{\hbox} \def\text{\hbox} \def\om{\omega} \def\Cal#1{{\bf #1}}
  \def\pcf{pcf} \defaultdefine{cf}{cf} \defaultdefine{reals}{{\Bbb R}}
  \defaultdefine{real}{{\Bbb R}} \def\restriction{{|}} \def\club{CLUB}
  \def\w{\omega} \def\exist{\exists} \def\se{{\germ se}} \def\bb{{\bf b}}
  \def\equivalence{\equiv} \let\lt< \let\gt>
\providecommand{\bysame}{\leavevmode\hbox to3em{\hrulefill}\thinspace}
\providecommand{\MR}{\relax\ifhmode\unskip\space\fi MR }
\providecommand{\MRhref}[2]{%
  \href{http://www.ams.org/mathscinet-getitem?mr=#1}{#2}
}
\providecommand{\href}[2]{#2}



\end{document}